\documentclass[english, makeidx]{amsart}

\usepackage{multicol} 
\usepackage{babel}
\usepackage{amstext}
\usepackage{amsmath}
\usepackage{amsfonts}
\usepackage{latexsym}
\usepackage{ifthen}
\usepackage{xypic}
\usepackage{amssymb}

\xyoption{all}
\newcommand{\Rien}[1]{}
\newcommand{\Formel}[1]{(\ref{#1})}

\newcommand{\B}{\mathbb B}

\newcommand{\Prop}[1]{Proposition~\ref{#1}}

\newcommand{\PN}{{\mathbb P}}

\newcommand{\lra}{\longrightarrow}
\newcommand{\KC}{{\mathbb C}}
\newcommand{\KZ}{{\mathbb Z}}
\newcommand{\KQ}{{\mathbb Q}}

\newcommand{\End}{{\rm End}}

\newcommand{\OH}{\mathbb H}
\newcommand{\KR}{\mathbb R}

\newtheorem{lemma1}[equation]{}

\newenvironment{example}{\begin{lemma1}{\bf Example.}\rm}{\end{lemma1}}

\newenvironment{proposition}{\begin{lemma1}{\bf Proposition.}}{\end{lemma1}}
\newenvironment{corollary}{\begin{lemma1}{\bf Corollary.}}{\end{lemma1}}

\newenvironment{rem}{\begin{lemma1}{\bf Remark.}\rm}{\end{lemma1}}

\begin{document}

\title{On abelian families and holomorphic normal projective connections}
\author[I. Radloff]{Ivo Radloff}
\address{I. Radloff - Mathematisches Institut - Universit\"at Bayreuth
  - D-95440 Bayreuth, Germany}
\email{ivo.radloff@uni-bayreuth.de}
\thanks{The author gratefully acknowledges support by the Schwerpunkt program {\em Globale Methoden in der komplexen Geometrie} of the Deutsche Forschungsgemeinschaft.}
\date{\today}
\maketitle

\section*{Introduction}
Our previous results from \cite{JRproj} raise the question whether the
list of complex
projective manifolds admitting a holomorphic normal projective
connection is the following (up to \'etale coverings):
\begin{enumerate}
 \item $\PN_m(\KC)$,
 \item smooth abelian families,
 \item manifolds with universal covering $\B_m(\KC)$.
\end{enumerate}
Here $\B_m(\KC)$ denotes the ball in $\KC^{m+1}$, the non
compact dual of $\PN_m(\KC)$ in the sense of hermitian symmetric
spaces. The second point inlcudes the flat case of an abelian
manifold.

Any compact Riemann surface admits a holomorphic normal projective
connection, this is the famous uniformization theorem. Kobayashi and
Ochiai showed that the list of projective surfaces with a holomorphic
normal projective connection consists of $\PN_2(\KC)$, abelian surfaces and
ball quotients (\cite{KO}). The above list was confirmed in
the case of projective threefolds in \cite{JRproj}.

Of particular interest are the manifolds with a holomorphic normal
projective connection of intermediate Kodaira dimension $0 < \kappa(M)
< m$. The type we
expect are locally symmetric spaces obtained as quotients of
  \[\KC^{m-n} \times \B_n(\KC)\]
by some special group of automorphisms as we will see in this
article. Concrete examples are given in the last section.

\section{Holomorphic normal projective connections}

Cartan's original definition of projective structures and connections
involve the language of principal bundles. We follow Kobayashi
and Ochiai (\cite{KO}). 

Let $M$ be some $m$--dimensional projective
manifold. Then $M$ carries a holomorphic normal projective connection if the (normalised) Atiyah class of the holomorphic cotangent bundle
has the form
  \begin{equation} \label{AtProj} 
a(\Omega_M^1) = \frac{c_1(K_M)}{m+1} \otimes id_{\Omega_M^1} + id_{\Omega_M^1}
  \otimes \frac{c_1(K_M)}{m+1} \in H^1(M, \Omega_M^1 \otimes T_M \otimes
  \Omega_M^1)
  \end{equation}
where we use the identities $\Omega_M^1 \otimes T_M \otimes
\Omega_M^1 \simeq
\End(\Omega_M^1) \otimes \Omega_M^1 \simeq \Omega_M^1 \otimes
\End(\Omega_M^1)$. The following Chern class identities hold, similar
to projective space:
\begin{equation} \label{chern}
  c_r(M) = \frac{1}{(m+1)^r}{m+1 \choose r}c_1^r(M) \quad \mbox{in }
H^r(M, \Omega_M^r)
 \end{equation}

\

It was shown in \cite{MM} that any holomorphic cocycle solution to
\Formel{AtProj} can be thought of as a $\KC$--bilinear map
  \[\Pi: T_M \times T_M \lra T_M\]
satisfying certain rules modelled after the Schwarzian derivative
\begin{enumerate}
  \item \[\Pi_{fX}Y = f\Pi_XY - \frac{1}{m+1}X(f)Y, \quad \mbox{f\"ur
    } f \in C^{\infty}(M)\]
  \item \[\Pi_X(fY) = f\Pi_XY + X(f)Y - \frac{1}{m+1}Y(f)X, \quad \mbox{f\"ur
    } f \in C^{\infty}(M)\]
  \item \[\Pi_XY - \Pi_YX = [X, Y]\]
\end{enumerate}
We shall not use this description in the sequel.

\

A manifold is said to carry a projective structure if there exists a
holomorphic projective atlas, i.e., an atlas with embeddings of the charts into some $\PN_m(\KC)$
such that the coordinate change is given by restrictions of
projective automorphisms. 

A manifold with a projective structure carries a (flat) projective connection.

\begin{example} \label{Exmpl}
  Projective space $\PN_m(\KC)$ carries a projective
  structure. Any manifold whose universal covering space admits an
embedding into $\PN_m(\KC)$ such that its fundamental group acts by
restrictions of projective transformations admits a projective
structure. In particular $\B_m(\KC)$, the non compact dual of $\PN_m(\KC)$, carries a
  projective structure. Any abelian manifold carries a projective structure.
Any Riemann surface carries a projective structure.
\end{example}

\

We call $\PN_m(\KC)$, ball quotients and {\'e}tale quotients of
abelian manifolds the {\em standard examples} of manifolds with a
projective structure. 

Let $M$ be a any {\em projective} manifold with a
holomorphic normal projective connection. If $K_M$ is not nef, then $M
\simeq \PN_m(\KC)$ (\cite{JRproj}). If $K_M$ is nef, then it is
expected that some multiple of $|K_M|$ is spanned defining a map
  \[f: M \lra N.\]
This is the famous abundance conjecture. Our results from
\cite{JRproj} suggest that for $M$ as above $f$ should be a
smooth abelian fibration, perhaps after some {\'e}tale covering.

\

We should mention that there are more examples if one drops the K\"ahler
condition. Twistor spaces over conformally flat Riemannian
fourfolds are complex threefolds with a projective structure. The only
K\"ahlerian twistor space is $\PN_3(\KC)$ by a result of
Hitchin. There are of course even more non--compact examples.

\section{Families of Abelian Varieties} \setcounter{equation}{0}
We are dealing with the following situation. Let $M$ be some
projective manifold of dimension $m$ admitting a smooth holomorphic
map
  \[f: M_m \lra N_n\]
onto some smooth projective manifold $N$ of positive dimension $n$. The fibers
are assumed to be $(m-n)$--dimensional abelian varieties.

We assume the existence of a smooth section. 

Projectivity of $M$ is not always necessary, but this is the case we
are interested in. The holomorphic one forms on
$M$ and $N$ give rise to the exact sequence
   \begin{equation} \label{reltang}
    0 \lra f^*\Omega_N^1 \stackrel{df}{\lra} \Omega_M^1 \lra
    \Omega^1_{M/N} \lra 0
   \end{equation}
where $\Omega^1_{M/N}$ is the sheaf of relative one forms. Here in
this case
\[E = f_*\Omega_{M/N}^1\]
is a vector bundle on $N$ of rank $m-n$ and $f^*E \simeq
\Omega_{M/N}^1$ via the canonical map $f^*f_*\Omega_{M/N}^1 \lra
\Omega_{M/N}^1$, as $\Omega^1_{M/N}$ is relatively spanned.

\begin{proposition} \label{fam}
  In the above situation, assume that $M$ admits a
  holomorphic normal projective connection. Then $N$ has a holomorphic normal
  projective connection, and
    \begin{equation} \label{AtE}
      a\big(E(-\frac{K_N}{n+1})\big) = 0 \quad
      \mbox{in } H^1(N, \End(E) \otimes
      \Omega_N^1),\end{equation}
  where $a$ denotes the Atiyah class of a vector bundle.
\end{proposition}

\begin{rem}
The formula is in terms of classes, we do not assume the
existence of a theta characteristic on $N$.
\end{rem}

\Prop{fam} will be proved below, we will first derive some consequences. The
trace of the Atiyah class gives the first Chern class, hence
  \begin{equation} \label{relK}
    c_1(E) = \frac{m-n}{n+1}c_1(K_N) \quad \mbox{in } H^1(N,
    \Omega_N^1).
  \end{equation}
Let as usual $K_{M/N}$ denote a divisor representing the determinant of 
$\Omega^1_{M/N}$. We have
  \[K_M = K_{M/N} + f^*K_N.\]
The divisor $K_{M/N}$ also represents $f^*\det E$ and \Formel{relK} gives
\begin{corollary}
  In the situation of the proposition the following identities hold in  $H^1(M, \Omega_M^1)$:
    \[c_1(K_{M/N}) = \frac{m-n}{n+1}c_1(f^*K_N) \quad \mbox{and} \quad
    c_1(K_M) = \frac{m+1}{n+1} c_1(f^*K_N)\]
  In particular, $c_1(K_M)$ and $c_1(f^*K_N)$ are proportional.
\end{corollary}

\begin{rem} 1.) The formulas hold in the case $m = n$. 

2.) If $c_1(K_M) = 0$ or $c_1(K_N) = 0$, then $c_1(K_M) = c_1(K_N) =
0$. By \Formel{chern} all Chern classes of $M$ and $N$ vanish. Then
$M$ and $N$ are abelian.

3.) If $K_M$ or $K_N$ is not nef, then $K_M$ and $K_N$ are not
nef. Then $M \simeq \PN_M(\KC)$ and $N \simeq \PN_n(\KC)$
(\cite{JRproj}). As $n>0$ by assumption, $m = n$ and
$f$ is an automorphism of projective space.

4.) In the case $\dim N = 1$ the formula reads
 \[2c_1(E) = 2c_1(f_*\Omega_{M/N}^1) = (m-1) c_1(K_N).\]
If $c_1(K_N) \not= 0$, i.e., in the non--abelian case, we get a family reaching
the Arakelov bound in the sense of Viehweg and Zuo in \cite{VZ}. 

Because of \Formel{AtE} and $c_1(K_N) \not= 0$, the bundle $E$ does
not admit a flat direct summand of positive rank. The result of \cite{VZ}
in this case seems to be that $M$ is isogeneous to some fiber product
 \[Z \times_Y \cdots \times_Y Z,\]
where $Z \lra Y$ is a modular family of abelian varieties with special
Hodge group built from quaternion algebras.
\end{rem}

\begin{proof}[Proof of \Prop{fam}]
  The arguments can essentially be found in \cite{JRproj}. By assumption we
  have a section $s: N \lra M$. 

Consider
  the pull back to $N$ by $s$ of \Formel{reltang}
   \begin{equation} \label{reltang2}
    0 \lra \Omega_N^1 \stackrel{s^*df}{\lra} s^*\Omega_M^1 \lra
    s^*\Omega^1_{M/N} \simeq E \lra 0.
   \end{equation}
  We have the map $ds: s^*\Omega_M^1 \lra \Omega^1_N$.   
 As $(ds)(s^*df) = d(f \circ s) = id_{\Omega_N^1}$, sequence \Formel{reltang2}
 splits holomorphically. 

\

The Atiyah class of
  $s^*\Omega_M^1$ is obtained from the Atiyah class of $\Omega^1_M$ by
  applying $ds$ to the last $\Omega_M^1$ factor in
  \Formel{AtProj}. What we get is 
    \begin{equation} \label{resat} 
    a(s^*\Omega_M^1) = \frac{s^*c_1(K_M)}{m+1} \otimes ds + id_{s^*\Omega_M^1} \otimes
    \frac{c_1(s^*K_M)}{m+1}
\mbox{ in } H^1(N, s^*\Omega_M^1 \otimes s^*T_M \otimes \Omega_N^1),
\end{equation} where
 we carefully distinguish between $s^*c_1(K_M) \in H^1(N,
 s^*\Omega_M^1)$ and the class $c_1(s^*K_M) = ds(c_1(K_M)) \in H^1(N, \Omega_N^1)$.

 The Atiyah class of a direct sum is the direct sum of the Atiyah
  classes. As the pull back of \Formel{reltang} splits
  holomorphically, we get the Atiyah classes of $\Omega_N^1$ and $E$ by
  projecting \Formel{resat} onto the corresponding summands. 

\

We begin with $E$. The class $c_1(K_M) \in H^1(M,
\Omega_M^1)$ is the pull back of some class in $H^1(N,
\Omega_N^1)$; it therefore vanishes under $H^1$ of $\Omega_M^1 \lra
\Omega^1_{M/N}$. This means the first summand in \Formel{resat}
vanishes if we project, while the second summand becomes
  \begin{equation} \label{AtEeins}
    id_E \otimes \frac{c_1(s^*K_M)}{m+1}  \;\; \mbox{in } H^1(N, E
  \otimes E^* \otimes
    \Omega_N^1).\end{equation}
This is the Atiyah class of $E$. The trace gives
  \[c_1(E) = rk(E)  \frac{c_1(s^*K_M)}{m+1} \in H^1(N,
  \Omega_N^1)\]
The determinant of
\Formel{reltang2} gives the following identities of classes in $H^1(N, \Omega_N^1)$:
  \begin{equation} \label{AdFor}
     c_1(K_N) = c_1(s^*K_M) - c_1(E) =
     \frac{m+1-(m-n)}{m+1}c_1(s^*K_M)
  \end{equation}
Now \Formel{AtE} follows from \Formel{AtEeins}.

\

We compute the Atiyah class of $\Omega_N^1$. We have to apply $ds$ to the
first factor of the first summand in \Formel{resat}. This gives
$\frac{c_1(s^*K_M)}{m+1}$. As the splitting maps compose to the identity we get
  \[\frac{c_1(s^*K_M)}{m+1} \otimes id_{\Omega_N^1} + id_{\Omega_N^1} \otimes
    \frac{c_1(s^*K_M)}{m+1}  \in H^1(N, \Omega_N^1 \otimes T_N \otimes
    \Omega_N^1)\]
This is the Atiyah class of $\Omega_N^1$. As we just saw in \Formel{AdFor}
  \[\frac{c_1(K_N)}{n+1} = \frac{c_1(s^*K_M)}{m+1}.\]
Replacing this in the above formula we see that $N$ has a holomorphic
normal projective connection. The proposition ist proved.
\end{proof}

\section{Examples} \setcounter{equation}{0}

The examples we give are well known PEL type Shimura
families. We follow the classical description of Shimura (\cite{Shim}) which makes the projective structure clearly visible.

The examples are quotients of
$\KC^{m-1} \times \B_1$. We identify $\B_1(\KC)$ and the upper half plane 
 \[\OH = \{\tau \in \KC \mid \Im m \tau > 0\}.\]
For any Fuchsian group $\Gamma \subset PGl_2(\KR)$ acting on
$\OH$ as a group of M\"obius transformations we denote by $\OH/\Gamma
= Y_{\Gamma}$ the corresponding quotient.

\subsection{Elliptic curves} A non--compact example and the case of a split algebra: let $B = M_{2\times
  2}(\KQ)$ and let $\Gamma$ be a torsion free congruence subgroup
of the group of positive units, i.e., of $Sl_2(\KZ)$. We let $\Gamma$ 
act on $\OH$ in the usual way. 

The elliptic curve $E_{\tau} = \KC/\Lambda_{\tau}$ where $\Lambda_{\tau} = \KZ \oplus
\KZ\tau$. If $\tau' =
\gamma(\tau) = \frac{a\tau + b}{c\tau + d}$ for some $\gamma \in
\Gamma$, then $E_{\tau} \simeq E_{\tau'}$. We will do the computation
below in the analogous case of false elliptic curves.  Consider the
subgroup $H_{\Gamma}$ of $Sl_3(\KC)$ of matrices
 \[\left(\begin{array}{c|c}
   \gamma & \begin{array}{c} 0 \\ 0 \end{array} \\ \hline
   \begin{array}{cc} m & n \end{array} & 1 \end{array}\right) \quad \mbox{with} \quad
 \gamma = \left(\begin{array}{cc} a & b \\ c & d \end{array}\right) \in
 \Gamma \mbox{ and } m, n \in \KZ.\]
The subgroup of maps where $\gamma = id$ is normal and isomorphic to 
$\Lambda \simeq \KZ \oplus \KZ$, and we obtain the
exact sequence 
  \begin{equation} \label{semdir}
    0 \lra \Lambda \lra H_{\Gamma} \simeq \Gamma \ltimes \KZ^{\oplus 2}  \lra
    \Gamma \lra 1.
  \end{equation}

Consider $\PN_2$ with homogeneous coordinates $\tau, u, z$.  Think of 
$\OH \times \KC$ with coordinates $\tau, z$ as an open subset
of the standard chart $u = 1$. The group $H_{\Gamma}$ induces a group
of projective automorphisms stabilizing $\OH \times \KC$. The above matrix acts as
  \[(\tau, z) \mapsto \left(\frac{a \tau + b}{c \tau + d}, \frac{z + m \tau
    + n}{c \tau + d}\right).\]
If we choose for $\Gamma$ some group such that $S_{\Gamma} = (\OH \times
\KC)/H_{\Gamma}$ becomes a smooth surface, then this surfaces has
a projective structure (example~\ref{Exmpl}). By
construction $S_{\Gamma}$ comes with an elliptic fibration over
$\OH/\Gamma$; the fiber over $[\tau]$ is isomorphic to the above
curve $E_{\tau}$.

Examples for $\Gamma$ are the well known congruence groups $\Gamma(N),
\Gamma_0(N), \Gamma_1(N)$ for certain level $N$. The surface $S_{\Gamma}$ is
not compact; any smooth compactification destroys the projective structure
(\cite{KO}).

\subsection{False elliptic Curves} \label{FEK}

\newcommand{\EO}{\mathfrak o}

A compact example and the case of a non split algebra: let $B$ be a total indefinite quaternion algebra
over $\KQ$. We fix an isomorphism $B \otimes \KR\simeq M_2(\KR)$. A
false elliptic curve is an abelian surface $F$ with
  \[\End(F) \otimes \KQ \simeq B.\]
The representing lattice is $\EO_B \simeq \KZ^4$, a maximal order in $B$. 
The orbit in $\KC^2$ of the vector $(\tau,
1)^t$, $\tau \in \OH$ under the matrices in $\EO_B$ is a complete
lattice $\Lambda_{\tau} \simeq \EO_B$. The quotient
$F_{\tau} = \KC^2/\Lambda_{\tau}$, for general $\tau$, is an example of a false elliptic
curve (\cite{Shim}).

Let $\Gamma$ be a torsion free congruence subgroup of the group of positive
units in $\EO_B$. Using $B \otimes \KR\simeq M_2(\KR)$ we 
consider $\Gamma$ as a group of matrices
and study its fixed point free action on $\OH$.

The quotient $Y_{\Gamma}$ is compact here (\cite{Shim}). If $\tau' =
\gamma(\tau) = \frac{a\tau + b}{c\tau + d}$ for some $\gamma \in \Gamma$, then 
 \[\Lambda_{\tau'} = \EO_B\left({\tau' \atop 1}\right) = \frac{1}{c
   \tau + d}\EO_B \left({a\tau+b \atop c\tau+d}\right) = \frac{1}{c
   \tau + d}\EO_B \gamma \left({\tau \atop 1}\right) = \frac{1}{c
   \tau + d}\Lambda_{\tau}\]
implying that $F_{\tau'}$ and $F_{\tau}$ are isomorphic.

\

As above we obtain a subgroup $H_B$ of $Sl_4(\KC)$ and an exact
sequence like \Formel{semdir} if we consider the
following matrices:
 \[\left(\begin{array}{c|c}
   \gamma_1 & 0 \\ \hline
   \gamma_2 & \begin{array}{cc}
                1 & 0 \\
                0 & 1
              \end{array} 
   \end{array}\right) 
\quad \mbox{with} \quad
\gamma_1 \in \Gamma, \gamma_2 \in \EO_B.\]

Consider $\PN_3$ with
homogeneous coordinates $\tau, u, z_1, z_2$.  We think of $\OH \times \KC^2$ with coordinates $\tau, z_1, z_2$ as an
open subset of the standard chart $u = 1$. The group $H_B$ induces a
group of projective automorphisms acting on $\OH \times \KC^2$. In
coordinates
  \[\gamma_k = \left(\begin{array}{cc} a_k & b_k \\
                                   c_k & d_k \end{array}\right), \quad
                               k = 1,2\]
act as
  \[(\tau, z_1, z_2) \mapsto \left(\frac{a_1 \tau + b_1}{c_1 \tau + d_1},
    \frac{z_1 + a_2\tau + b_2}{c_1 \tau + d_1}, \frac{z_2 + c_2\tau + d_2}{c_1 \tau + d_1}\right).\]
We see that $H_B$ acts on the first factor by $\gamma_1 \in \Gamma$. The quotient
  \[M_{B, \Gamma} = (\OH \times \KC^2)/H_B\]
carries the structure of a smooth projective threefold and has a
projective structure (example~\ref{Exmpl}). The fiber of 
$M_{B, \Gamma} \lra Y_{\Gamma}$  over $[\tau]$ is isomorphic to the above false elliptic
curve $F_{\tau}$. 

\

The latter construction also works if we consider the split algebra $M_{2
  \times 2}(\KQ)$ instead of a quaternion algebra. The quotient in
this case, taking for $\Gamma$ the same group as above, is just the
fiber product $S_{\Gamma} \times_{Y_{\Gamma}} S_{\Gamma}$ from
above. Likewise, if we consider
\[\left(\begin{array}{c|c|c}
   \gamma_1 & 0 & 0 \\ \hline
   \gamma_2 & \begin{array}{cc}
                1 & 0 \\
                0 & 1
              \end{array} & 0 \\ \hline
    \gamma_3 & 0 & \begin{array}{cc}
                1 & 0 \\
                0 & 1
              \end{array}
   \end{array}\right),
\quad \quad
 \gamma_1 \in \Gamma, \gamma_2, \gamma_3 \in \EO_B.\]
acting on $\OH \times \KC^2 \times \KC^2$ we get $M_{B, \Gamma}
\times_{Y_\Gamma} M_{B, \Gamma}$, a compact example of a smooth
projective manifold with a projective structure.

\end{document}